\documentclass[11pt]{article}
\usepackage{amsmath}
\usepackage{amssymb}
\usepackage{url}

\newcommand{\id}[1]{\textnormal{\texttt{Id}}_{#1}}
\newcommand{\UU}{\mathcal{U}}
\newcommand{\eq}{\texttt{Eq}}

\begin{document}

\title{Voevodsky's Univalence Axiom in homotopy type theory}

\author{Steve Awodey, \'Alvaro Pelayo and Michael A. Warren}

\date{}

\maketitle

\begin{abstract}
In this short note we give a glimpse of homotopy type theory, a new field of mathematics
at the intersection of algebraic topology and mathematical logic, and we explain Vladimir Voevodsky's univalent 
interpretation of it.   This interpretation has given rise to the \emph{univalent foundations program}, which is the topic
of the current special year at the Institute for Advanced Study.
 \end{abstract}

 The Institute for Advanced Study in Princeton is hosting
a special program during  the academic year 2012\--2013 on a new research theme that is based on recently discovered connections  
between homotopy theory, a branch of algebraic topology, and type theory, a branch of mathematical logic and theoretical computer science. 
In this brief paper our goal is to take a glance at these developments. For those readers who would like to learn more about them, we recommend a number of references throughout.

{\bf Type theory} was invented by Bertrand Russell \cite{Russell:1908}, but it was first developed as a rigorous formal system by
Alonzo Church \cite{Church:1933cl,Church:1940tu,Church:1941tc}.  It now has numerous applications in computer science, especially in the theory
of programming languages \cite{Pierce:2002tp}.   Per Martin-L\"{o}f
\cite{MartinLof:1998tw,MartinLof:1975tb,MartinLof:1982bn,MartinLof:1984tr}, among others,
developed a generalization of Church's system which is now usually called 
dependent, constructive, or simply {\bf Martin\--L\"of type theory}; this is the system that we consider here. It was originally
intended as a rigorous framework for constructive mathematics. 

In type theory objects are classified using a primitive notion of \emph{type}, similar to the data-types used in programming languages.  And as in programming languages, these elaborately structured types can be used to express detailed specifications of the objects classified, giving rise to principles of reasoning about them.  To take a simple example, the objects of a product type $A\times B$ are known to be of the form $\langle a, b\rangle$, and so one automatically knows how to form them and how to decompose them.
This aspect of type theory has led to its extensive use in verifying the correctness of computer programs.  Type theories also form the basis of modern computer proof assistants, which are used for formalizing mathematics and verifying the correctness of formalized proofs.  For example, the powerful Coq proof assistant \cite{coq} has recently been used to formalize and verify the correctness of the proof of the celebrated Feit-Thompson Odd-Order theorem \cite{gonthier}.

One problem with understanding type theory from a mathematical point of view, however, has always been that the basic concept of \emph{type} is unlike that of \emph{set} in ways that have been hard to make precise. This difficulty has now been solved by the idea of regarding types, not as strange sets (perhaps constructed without using classical logic), but as spaces, regarded from the perspective of homotopy theory

In {\bf homotopy theory} one is concerned with spaces and continuous mappings between them, 
up to homotopy; a \emph{homotopy} between a pair of continuous maps $f \colon X	\to Y$
and  $g \colon X	\to Y$ is 
a continuous map $H \colon X \times [0, 1]	\to Y$ satisfying
$H(x, 0) = f (x)$  and $H(x, 1) = g(x)$. The homotopy $H$ may be thought of as a ``continuous deformation" of $f$ into $g$. The spaces $X$ and $Y$ are said to be \emph{homotopy equivalent}, $X\simeq Y$, if there are continuous maps going back and forth, the composites of which are homotopical to the respective identity mappings, i.e.\ if they are isomorphic ``up to homotopy".  Homotopy equivalent spaces have the same algebraic invariants (e.g.\ homology, or the fundamental group), and are said to have the same \emph{homotopy type}.

{\bf Homotopy type theory} is a new field of mathematics which interprets type theory from a homotopical perspective.
In homotopy type theory, one regards the types as spaces, or homotopy types, and the logical constructions (such as the product $A\times B$) as homotopy-invariant constructions on spaces.  
In this way, one is able to manipulate spaces directly, without first having to develop point-set topology or even define the real numbers.
Homotopy type theory is connected to several topics of interest in modern algebraic topology, such as $\infty$-groupoids and Quillen model structures (see \cite{PeWa2012});  we will only mention one simple example below, namely the homotopy groups of spheres.

To briefly explain the homotopical perspective of types, consider the basic concept of type theory, namely that
the \emph{term} $a$ is of \emph{type} $A$, which is written:
$$
  a:A.
$$
This expression is traditionally thought of as akin to ``$a$ is an
element of the set $A$."  However, in homotopy type theory we think of
it instead as ``$a$ is a point of the space $A$."  Similarly, every term $f : A \to B$ is regarded as a continuous function from the space $A$ to the space $B$. 

This perspective clarifies features of type theory which were puzzling from the perspective of types as sets; for instance, that one can have non-trivial types $X$ such that $(X\to X) \cong X + 1$.  But the key new idea of the homotopy interpretation is that the logical notion of identity $a = b$ of two objects $a, b: A$ of the same type $A$ can be understood as the existence of a path $p : a \leadsto b$ from point $a$ to point $b$ in the space $A$.  This also means that two functions $f, g: A\to B$ are identical just in case they are homotopic, since a homotopy is just a family of paths $p_x: f(x) \leadsto g(x)$ in $B$, one for each $x:A$.  In type theory, for every type $A$ there is a (formerly somewhat mysterious) type $\id{A}$ of identities between objects of $A$; in homotopy type theory, this is just the \emph{path space} $A^I$ of all continuous maps $I\to A$ from the unit interval. (See  \cite{Awodey:2009bz,Aw2010,PeWa2012}.) 

At around the same time that Awodey and Warren advanced the idea of homotopy type theory, Voevodsky showed how to model type theory using Kan simplicial sets, a familiar setting for classical homotopy theory, thus arriving independently at essentially the same idea around 2005.  Both were inspired by the prior work of Hofmann and Streicher, who had constructed a model of type theory using groupoids \cite{HS}.

Voevodsky moreover recognized that this simplicial interpretation satisfies a further crucial property, which he termed \emph{univalence},  and which is not usually assumed in type theory.  Adding univalence  to type theory in the form of a new axiom has far-reaching consequences, many of which are natural, simplifying and compelling.  The {\bf Univalence Axiom} thus further strengthens the homotopical view of type theory, since it holds in the simplicial model, but fails in the view of types as sets. 

The basic idea of the Univalence Axiom can be explained as follows.  In type theory, one can have a universe $\UU$, the terms of which are themselves types, $A : \UU$, etc.  Of course, we do not have $\UU:\UU$, so only some types are terms of $\UU$ -- call these the \emph{small} types.  Like any type, $\UU$ has an identity type $\id{\UU}$, which expresses the identity relation $A = B$ among small types.  Thinking of  types as spaces, $\UU$ is a space, the points of which are spaces; to understand its identity type, we must ask, what is a path $p : A \leadsto B$ between spaces in $\UU$?  The Univalence Axiom says that such paths correspond to homotopy equivalences $A\simeq B$, as explained above (the actual notion of equivalence required is slightly different).  A bit more precisely, given any (small) types $A$ and $B$, in addition to the type $\id{\UU}(A,B)$ of identities between $A$ and $B$ there is the type $\eq(A,B)$ of equivalences from $A$ to $B$.  Since the identity map on any object is an equivalence, there is a canonical map,
$$\id{\UU}(A,B)\to\eq(A,B).$$
The Univalence Axiom states that this map is itself an equivalence.  At the risk of oversimplifying, we can state this succinctly as follows:

\begin{description}
\item[Univalence Axiom]  $(A = B)\ \simeq\ (A\simeq B)$.
\end{description}
In other words, identity is equivalent to equivalence. 

From the homotopical point of view, this says that the universe $\UU$ is something like a classifying space for (small) homotopy types, which is a practical and natural assumption.  From the  logical point of view, however, it is 
revolutionary: it says that isomorphic things can be identified!  Mathematicians are of course used to identifying isomorphic structures in practice, but they generally do so with a wink, knowing that the identification is not ``officially" justified by foundations.  But in this new foundational scheme, not only are such structures formally identified, but the different ways in which such  identifications may be made themselves form a structure that one can (and should!) take into account.

Part of the appeal of homotopy type theory with the Univalence Axiom is the many interesting connections it reveals between logic and homotopy.  Another remarkable aspect is that it can be carried out in a {\bf computer proof assistant}, since type theory exhibits such good computational properties (see \cite{Simpson:2004bt,Hales:2008ud} on the use of computer proof assistants in general).  In practical terms, this means that it is possible to use the powerful, currently available proof assistants based on type theory, like the Coq system, to develop mathematics involving homotopy theory, to verify the correctness of proofs, and even to provide some degree of automation of proofs.  

To give just one example, in homotopy type theory one can directly define the $n$-dimensional sphere $S^n$ as a type, with its associated principles of reasoning.  Moreover, for any type $A$ one can define the homotopy groups $\pi_n(A)$, again in a very direct way in terms of the identity type $\id{A}$ explained above.   One can then reason directly in type theory, using the principles associated with these constructions, and prove for example that $\pi_n(S^n) = \mathbb{Z}$ for $n\geq 1$ (as has recently been done by G.~Brunerie and D.~Licata at the Institute for Advanced Study, using the Univalence Axiom in an essential way).  Finally, the proof can be  formalized in a proof assistant and verified by a computer.  In this way, one not only has new methods of proof in classical homotopy theory, but indeed ones which provide associated computational tools.

Voevodsky has christened this combination of homotopy type theory with the Univalence Axiom, implemented on a computer proof assistant, the {\bf Univalent Foundations} program.  It can be regarded as a new foundation for mathematics in general, not just for homotopy theory, as Voevodsky has shown by developing an extensive code library of formalized mathematics in this setting.  Moreover, he is promoting more interaction between pure mathematicians and the developers of such proof assistants, as is occurring in the special year on Univalent Foundations at the Institute for Advanced Study.  

For those interested in contributing to this new kind of mathematics, it may be encouraging to know that there are many interesting open questions.  The most pressing of them is perhaps the ``constructivity'' of the Univalence Axiom itself, conjectured by Voevodsky in \cite{Vo2012}.  
It concerns the effect of adding the Univalence Axiom on the computational behavior of the system of type theory, and thus on the existing proof assistants.  Another major direction, of course, is the further formalization of classical results and current mathematical research in the univalent setting.  We expect that it will eventually be possible to formalize large amounts of modern mathematics in this setting, and that doing so will give rise to both theoretical insights and good numerical algorithms (extracted from code in a proof assistant).  

In this direction, together with Voevodsky, the last two authors are
working on an approach to the theory of integrable systems (using the
new notion of $p$-adic integrable system as a test case) in the
univalent setting. A preliminary treatment is the construction of the
$p$\--adic numbers is given in \cite{PeVoWa2012}.  One of Voevodsky's
goals (as we understand it) is that in a not too distant future,
mathematicians will be able to verify the correctness of their own
papers by working within the system of univalent foundations
formalized in a proof assistant, and that doing so will become
natural even for pure mathematicians (the same way that most
mathematicians now typeset their own papers in \TeX).  We believe
that this aspect of the univalent foundations program distinguishes
it from other approaches to foundations, by providing a practical
utility for the working mathematician.

Our goal in this announcement has been to give a brief and intentionally superficial glimpse of two closely related, recent developments: homotopy type theory and Voevodsky's univalent foundations program. At present, these subjects are still developing quite rapidly, and the current literature is (with few exceptions) highly specialized and, unfortunately, largely inaccessible to those without prior knowledge of homotopy theory and logic. One exception is the survey article  \cite{PeWa2012}, which goes into much greater depth than the present article, while still being intended for a general mathematical readership; it also contains an introduction to the use of the Coq proof assistant in the univalent setting. See also \cite{Aw2010,HoTT,Vo2012}.

\vspace{5mm}
\emph{Acknowledgements}. 
We thank the Institute for Advanced Study for the
excellent resources which have been made available to the authors during the
preparation of this article. We thank  Thierry Coquand, Dan Grayson, and 
Vladimir Voevodsky for useful discussions on the topic of this paper, and the referees
for helpful suggestions. Awodey is partly supported by National Science 
Foundation Grant DMS-1001191 and Air Force OSR Grant 11NL035. Pelayo is partly supported by NSF 
CAREER Grant DMS-1055897, NSF Grant DMS-0635607, and Spain Ministry of Science Grant 
Sev-2011-0087. Warren is supported by the Oswald Veblen Fund. 

\providecommand{\bysame}{\leavevmode\hbox to3em{\hrulefill}\thinspace}

\noindent
\\
{\bf Steve Awodey}\\
 Carnegie Mellon University\\
 Department of Philosophy\\
 Pittsburgh, PA  15213 USA
 \\
  \\
  Institute for Advanced Study\\
 School of Mathematics, Einstein Drive\\ 
 Princeton,  NJ 08540 USA \\
E-mail addresses: awodey@cmu.edu, awodey@math.ias.edu
\\
\\
\\
{\bf \'Alvaro Pelayo}\\
Washington University\\  
   Department of Mathematics\\
   One Brookings Drive, Campus Box 1146\\
   St Louis, MO 63130 USA
   \\
   \\
   Institute for Advanced Study\\
 School of Mathematics, Einstein Drive\\
 Princeton,  NJ 08540 USA\\
 E-mail addresses: apelayo@math.wustl.edu, apelayo@math.ias.edu
 \\
 \\
 \\
{\bf Michael A. Warren}\\
Institute for Advanced Study\\
 School of Mathematics, Einstein Drive\\
 Princeton, NJ 08540 USA.\\
 E-mail address: mwarren@math.ias.edu\\

\end{document}